%% file: article.tex
\documentclass[12pt]{amsart}

\usepackage[utf8]{inputenc}
\usepackage[english]{babel}
\usepackage{amsmath,amssymb} 
\usepackage[usenames, dvipsnames]{xcolor}

\usepackage{mystyle}
\usepackage{myvisual}
\usepackage{notations}

\setulcolor{gray} 

\title{An Epigraphical Approach to the Representer Theorem}
\author{Vincent Duval  \\
	INRIA \& CEREMADE, Universit\'e Paris-Dauphine, PSL  \\
	}

\date{\today}
\begin{document}

\maketitle

\begin{abstract}
Describing the solutions of inverse problems arising in signal or image processing is an important issue both for theoretical and numerical purposes. We propose a principle which describes the solutions to convex variational problems involving a finite number of measurements. We discuss its optimality on various problems concerning the recovery of Radon measures.
\end{abstract}
\vspace{0.5cm}

The recovery an unknown vector from a finite number of linear measurements is an inverse problem which frequently arises  
in signal processing or machine learning. Convex variational approaches provide a flexible framework to address that task, typically by solving a problem of the form
\begin{align}\label{faces:eq:minconv}
  \min_{u\in \vecgal} \reg(u) + \fido(\Phi u), \tag{$\Pp$}
\end{align}
where $\vecgal$ is a real vector space (\eg{} the space of signals), $\Phi:\vecgal \rightarrow \RR^\nbmes$ is linear, and $\reg:\vecgal \to \RR\cup\{+\infty\}$ and  $\fido:\RR^\nbmes \rightarrow \RR\cup \{+\infty\}$ are two convex functions.

The solutions to~\eqref{faces:eq:minconv} are in general not unique, and it is of crucial importance to understand the structure of the solution set $\sol$. As noted in~\cite{chandrasekaran2012convex}, the regularization term $\reg$ determines the building blocks, or \emph{atoms}, which constitute, by convex combinations, the (or some) elements of $\sol$. Those atoms are the extreme points of the level sets of $\reg$. Besides its theoretical interest, that description makes a strong case for greedy optimization algorithms such as the conditional gradient (a.k.a. Frank-Wolfe)~\cite{bredies_inverse_2013,boyd-adcg2015}, especially in the case where $\vecgal$ is infinite-dimensional: even though it is not possible to encode a generic vector of $\vecgal$ on a computer, the \emph{atoms} induced by $\reg$ might be simple enough to be handled (possibly with some approximation) numerically.

The goal of the present note is to state this representation principle as precisely as possible, that is, the description of solutions of~\eqref{faces:eq:minconv} as a convex combination of \emph{atoms}, while emphasizing the geometric essence of the property.

\section{Related works and main theorem}
\subsection{Representer theorems.}
While representation results describing the sparsity (resp. rank) of some solutions of linear (resp. semi-definite) programs are well-known (see \eg{} \cite{matousek2007understanding,barvinok_problems_1995}), our main focus is on infinite dimensional problems, for instance the \emph{Basis Pursuit for measures}~\cite{de_castro_exact_2012,candes2014towards},
\begin{equation}\label{eq::beurling}
  \min_{m\in \bRad} \tvmes{\mes}  \quad \mbox{s.t.}\quad  \Phi \mes=y 
\end{equation}
where $y\in\RR^\nbmes$, $\domX$ is a compact subset of $\RR^d$ or the $d$-dimensional torus $\TT^d$, $\bRad$ is the space of Radon measures on $\domX$ and $\tvmes{\mes}$ denotes the total variation of $\mes$. The observation is $\Phi \mes=\left(\int_\domX \varphi_i(x)\D\mes(x)\right)_{1\leq i\leq \nbmes}$ where $\{\varphi_i\}_{1\leq i\leq \nbmes}$ is a family of continuous functions on $\domX$. It is known~\cite{Zuhovickii1948} that there is a solution which can be represented as  
\begin{align}
  \label{eq:sumdirac}
  \mes=\sum_{i=1}^r a_i \delta_{x_i} \quad \mbox{where $r\leq \nbmes$, $\{(x_i, a_i)\}_{i=1}^r\subseteq \Omega\times \RR$.}
\end{align} More precisely, Fisher and Jerome~\cite{fisher_spline_1975} (see also~\cite{unser2017splines,flinth_exact_2019}) have proved that \emph{the extreme points of the solution set $\sol$ are of the form~\eqref{eq:sumdirac}}. Since, in a locally convex Hausdorff space, each nonempty compact convex set is the closed convex hull of its extreme points, knowing the solutions of the form~\eqref{eq:sumdirac} allows to recover the whole set of solutions.

A related result was proved by L.C.~Dubins~\cite{dubins_extreme_1962} and V.~Klee~\cite{klee_theorem_1963} in the study of convex sets: the extreme points of the intersection of a linearly closed and bounded convex set $C$ with an affine space of codimension $\nbmes$ are a convex combination of at most $\nbmes+1$ extreme points of $C$.
In a variational problem with an equality constraint such as~\eqref{eq::beurling} (\ie{} $\fido(w)=\chi_{\{y\}}(w)\eqdef 0$ if $w=y$, $+\infty$ otherwise), one may regard the solution set as the intersection of the level set $\lst{\reg}{\min \Pp}$ with $\Phi^{(-1)}\{y\}$. 
However, compared to~\eqref{eq:sumdirac}, applying the Dubins-Klee theorem yields one atom too many, and in~\cite{boyer_representer_2019} (together with coauthors) we have provided an argument which reduces that number to $\nbmes$ by exploiting the structure of convex optimization problems (see also~\cite{bredies_sparsity_2018}).

The present note describes an alternative argument which handles the case of general convex fidelity terms much more precisely than in~\cite{boyer_representer_2019} and makes the statement more symmetric. Moreover, we discuss the optimality of the theorem in the particular case of Radon measure recovery.

\subsection{Main theorem.}
Our main result relates the dimensions of the faces in the level sets of the regularizer and the fidelity term,
\begin{align*}
  \lst{\reg}{\reg(p)} &\eqdef \enscond{u\in \vecgal}{\reg(u)\leq \reg(p)},\\ 
  \lst{\fido}{\fido(\Phi p)} &\eqdef \enscond{w \in \RR^\nbmes}{\fido(w)\leq \fido(\Phi p)}.
\end{align*}
The notion of face used below is described in~\Cref{sec:faceconv}. Note that, possibly changing $f$ and $\nbmes$, there is no loss of generality in assuming that $\Phi$ is surjective.
\begin{theorem}
\label{thm:first}
  Let $\reg:\vecgal \to \RR\cup\{+\infty\}$, $\fido:\RR^\nbmes \rightarrow \RR\cup \{+\infty\}$ be two convex functions, and let $\Phi:\vecgal \rightarrow \RR^\nbmes$ be linear surjective. Assume that $p \in \sol$, $\reg(p)+\fido(\Phi p)<+\infty$, and that $\lst{\reg}{\reg(p)}$ is linearly closed and contains no line.  
  
  If $p$ belongs to a face of $\sol$ with dimension $j<+\infty$, then it belongs to a face of $\lst{\reg}{\reg(p)}$ with dimension at most $\dimreg$, where
  \begin{align}
    k\eqdef \begin{cases}
      \nbmes -\dimfidlst+j-1 & \mbox{if $\left(\reg(p)>\inf \reg\right)$ or $\left(\fido(\Phi p)>\inf \fido\right)$},\\
      \nbmes -\dimfidlst+j & \mbox{otherwise.}
    \end{cases}
  \end{align}
  and $\ell$ is the dimension of the minimal face of $\Phi p$ in $\lst{\fido}{\fido(\Phi p)}$.  If, moreover, $p$ satisfies the double obliqueness condition described in~\Cref{faces:def:oblique}, the number $k$ can be reduced to $\nbmes -\dimfidlst+j-2$. 
 \end{theorem} 
We deduce the following representation for $p$.
\begin{corollary}\label{coro:repres}
Under the assumptions of \Cref{thm:first}, $p$ can be written as a convex combination of (at most)
    $\dimreg+1$ extreme points of $\lst{\reg}{\reg(p)}$, or $\dimreg$ points of $\lst{\reg}{\reg(p)}$, each an extreme point or a point in an extreme ray.
\end{corollary}

\begin{rem}
  In the case of an equality constraint ($\fido=\chi_{\{y\}}$) or a strictly convex function $\fido$, we have $\dimfidlst=0$. Hence, $p$ belongs to a face of $\lst{\reg}{\reg(p)}$ with dimension (at most) $\nbmes +j-1$ if it belongs to a face of $\sol$ with dimension $j$ and $\reg(p)<\inf \reg$. This is coherent with~\cite{boyer_representer_2019}. On the other hand, with polyhedral fidelity terms $\fido$ such as the $\ell^1$ or $\ell^\infty$ norm, non-trivial values of $\ell$ should be taken into account.
\end{rem}
\begin{rem}
 Furthermore, for $j=0$, \ie{} $p$ extreme point of $\sol$, we obtain that $p$ belongs to a face of dimension (at most) $\nbmes-1$ provided $\reg(p)>\inf \reg$. Hence \Cref{coro:repres} recovers the Fisher-Jerome theorem.
\end{rem}

\section{Proof of the main result}

\subsection{Reminder on the faces of convex sets.}\label{sec:faceconv} We first recall a few definitions and basic properties. 
Given two points $x$ and $y$ in $\vecgal$, we define the
closed interval (or line segment) joining $x$ to $y$ as $\ci{x}{y}\eqdef\enscond{tx+(1-t)y}{0\leq t\leq 1}$, and the open interval joining $x$ to $y$ as {$\oi{x}{y}\eqdef \ci{x}{y} \setminus \{x,y\}$}. 
A line (resp. an open half line) is a set of the form $a+\RR v$ (resp. $\enscond{a+tv}{t>0}$) where $a,v\in \vecgal$ and $v\neq 0$. 
In the following, $\cvx\subseteq \vecgal$ denotes a \emph{convex set}, \ie{} for any $x,y\in \cvx$, the segment $\ci{x}{y}$ lies in $\cvx$.

The set $\cvx$ is \emph{linearly closed} (resp. linearly bounded) if its intersection with any line is a closed (resp. bounded) subset of that line.

 We say that a point $u\in \cvx$ belongs to the \emph{relative algebraic interior} (or \emph{intrinsic core}) of $\cvx$  if 
\begin{align}
  \forall v\in \Aff \cvx,\quad \exists \epsilon>0,\ \forall \lambda\in \oi{-\epsilon}{\epsilon},\quad u+\lambda(v-u)\in \cvx,
\end{align}
where $\Aff \cvx$ denotes the affine hull of $\cvx$. Equivalently (see~\cite[Prop. 2.3]{klee_extremal_1957}), $u$ is in the relative algebraic interior of $\cvx$ if and only if 
\begin{align}\label{faces:eq:internal}
  \forall v\in \cvx\setminus \{u\},\ \exists z\in \cvx,\ u\in \oi{v}{z}.
\end{align}
We say that $\cvx$ is \emph{internal} if it is equal to its relative algebraic interior.

 A point $x\in \cvx$ is an \emph{extreme point} of $\cvx$ if there is no open interval in $\cvx$ containing $x$, or equivalently if $\cvx\setminus\{x\}$ is convex.
An \emph{extreme ray} $\rho$ of $\cvx$ is a half-line contained in $\cvx$ such that any open interval $I$ which intersects $\rho$ must satisfy $I\subseteq \rho$.
More generally, a subset $F$ of $\cvx$ is said to be a \emph{face} of $\cvx$ if $F$ is convex and \revision{for all  open interval $I\subseteq \cvx$ which intersects $F$,} $I\subset F$. An alternative definition of an extreme point is ``a point $x$ such that $\{x\}$ is a face of $\cvx$''. Similarly, extreme rays may be defined as the half-lines which are a face of $\cvx$.
The dimension of a face, $\dim F$, is defined as the dimension of its affine hull $\Aff(F)$.

A canonical choice of face is given by the notion of \emph{elementary face}. Given $x\in\cvx$, we define $\face{x}{\cvx}$ as the intersection of all the faces of $\cvx$ which contain $x$. It is also a face, hence it is the \emph{minimal face} of $\cvx$ (for the inclusion) which contains $x$. We call such sets the elementary faces of $\cvx$. It turns out that $\face{x}{\cvx}$ is equal to the largest internal subset of $\cvx$ which contains $x$ (see~\cite[Th. 2.1]{dubins_extreme_1962}),  hence it is the union of $\{x\}$ and all the open intervals of $\cvx$ which contain $x$. Moreover, the elementary faces yield \emph{a partition of $\cvx$} with $y\in \face{x}{C}$ if and only if $\face{x}{\cvx}=\face{y}{C}$.

The behavior of the elementary faces when performing several operations on convex sets is described below.

\emph{Intersection.}
Since the elementary face $\face{x}{\cvx}$ is the union of $\{x\}$ and all the open intervals of a convex set which contain $x$, one may check that if $\cvx_1$ and $\cvx_1$ are two convex sets,
\begin{align}
  \label{eq:faceintercvx}
  \face{x}{\cvx_1\cap \cvx_2}&=\face{x}{\cvx_1}\cap \face{x}{\cvx_2}.
\end{align}
Moreover, if $W_{1,2}$, $W_1$, $W_2$  respectively denote the affine hulls of those faces, they consist in the collection of lines through $x$ which respectively intersect $\cvx_1\cap \cvx_2$, $\cvx_1$, $\cvx_2$ through an open interval. As a consequence,\begin{align}
  \label{eq:spanintercvx}
  W_{1,2}&=W_1\cap W_2.
\end{align}

\emph{Cartesian product.}
If $\cvx_1$, $\cvx_2$ are convex subsets of the vector spaces $\vecgal_1$, $\vecgal_2$, it is possible to check that $\face{x_1}{\cvx_1}\times \face{x_2}{\cvx_2}$ is both a face of $\cvx_1\times \cvx_2$ and an internal set. As a result,
\begin{align}
  \label{eq:faceprodcvx}
  \face{(x_1,x_2)}{\cvx_1\times \cvx_2}&=\face{x_1}{\cvx_1}\times \face{x_2}{\cvx_2}.
\end{align}
Moreover, if $W_{1,2}$, $W_1$, $W_2$  respectively denote the affine hulls of the above-mentioned faces, it holds
\begin{align}
  \label{eq:spanprodcvx}
  W_{1,2}&=W_1\times W_2.
\end{align}

\emph{Affine map.}
If $\psi:\vecgal_1\rightarrow \vecgal_2$ is an affine bijective map, it preserves the elementary faces:
\begin{align}
  \label{eq:facebijcvx}
  \face{\psi(x)}{\psi(\cvx)}&=\psi\left(\face{x}{\cvx}\right).
\end{align}
If $W_1$ (resp. $W_2$) denotes the affine hull $\face{x}{\cvx}$ (resp. $\face{\psi(x)}{\psi(\cvx)}$),
\begin{align}
  W_2&=\psi(W_1). 
\end{align}

\subsection{Epigraphical reformulation}\label{sec:epireform} Now, we proceed with the proof of~\Cref{thm:first}. We consider $p$ as in~\Cref{thm:first}, and we write $\fidp\eqdef \fido\circ\Phi$, $\minval\eqdef \reg(p)+\fidp(p)$. Instead of directly studying $\sol$ and the level sets $\lst{\reg}{\reg(p)}$, $\lst{\fidp}{\fidp(p)}$, we consider $\solh \eqdef \enscond{(u,\reg(u))}{u\in \sol}$, as well as the epigraph of $\reg$ and the hypograph of $\minval-\fidp$,
\begin{align*}
  \epir&\eqdef  \epi(\reg)\eqdef \enscond{(u,r)\in \vecgal\times \RR}{\reg(u)\leq r},\\
  \hypof&\eqdef \hypo(\minval-\fidp)\eqdef \enscond{(u,r)\in \vecgal\times \RR}{\minval-\fidp(u)\geq r}.
\end{align*}
We note that
$\epir$ and $\hypof$ are convex, with \todetail{$\solh=\epir\cap \hypof$} (see \cref{fig:epigraphs} for an illustration).
\begin{detail}
  Indeed, if $(u,r)\in \epir\cap \hypof$, then $\reg(u)\leq r\leq \minval-\fidp(u)$, hence
  \begin{align*}
    \reg(u)+\fidp(u)\leq r+\fidp(u)\leq \minval.
  \end{align*}
  Since $\minval=\inf\eqref{faces:eq:minconv}$, the left-hand side is bounded below by $\minval$, hence $u\in \sol$ and $r=\reg(u)$. This proves that $\epir\cap \hypof\subseteq \solh$. 
  The converse inclusion is straightforward.
\end{detail}

  Using~\eqref{eq:faceintercvx}, we get  \begin{align}
  \label{eq:faces}
  \face{\vecreg{p}}{\epir\cap \hypof}=\face{\vecreg{p}}{\epir}\cap \face{\vecreg{p}}{\hypof}.
\end{align}
  To understand the dimension of those faces, we need to consider their affine hulls. Let $\hS$, $\hE$ and $\hH$ be the respective affine hulls of $\face{\vecreg{p}}{\epir\cap \hypof}$, $\face{\vecreg{p}}{\epir}$ and $\face{\vecreg{p}}{\hypof}$.
 From~\eqref{eq:spanintercvx}, we have $\hS= \hE\cap \hH$. 
  Up to a translation of the origin in $\vecgal\times \RR$, we assume without loss of generality that $\vecreg{p}=\vecepi{0}{0}$, so that \emph{all the above-mentioned affine hulls are now linear hulls}. By classical results in linear algebra, 
  $\dim \hE = \dim(\hE\cap \hH) +  \codim_{\hE}(\hE\cap \hH)$ and $\codim_{\hE}(\hE\cap \hH)=\codim_{\hE+\hH}(\hH)= \codim_{\vecgal\times \RR} \hH - \codim_{\vecgal\times \RR}(\hE+\hH)$.
Combining those equalities, we get 
\begin{align}\label{faces:eq:dimensions}
  \dim \hE = \codim_{\vecgal\times\RR} \hH + \dim(\hS) - \codim_{\vecgal\times\RR}(\hE+\hH).
\end{align}

\begin{figure}
  \centering
  \tdplotsetmaincoords{85}{22}
    \input{fig-epigraphs}
    \caption{The solution set $\sol$ is equivalent, up to an affine isomorphism, to the set $\solh$ (see \cref{faces:sec:linear}).}
  \label{fig:epigraphs} 
  \vspace{0.5cm}
\end{figure}

  \subsection{Faces in epigraph and in level sets}\label{sec:facepi} To relate the faces of $\solh$ and $\sol$, we note that $\solh$ is the image of $\sol$ by some injective linear map. 
  \begin{lemma}\label{faces:sec:linear}
  There is a linear  map $\mapS:\Span(\sol)\rightarrow \RR$ such that $\reg$ coincides with $\mapS$ on $\sol$.
  Moreover, the map $\hmapS:\Span(\sol)\rightarrow \Span(\solh)$ defined by $\hmapS(v)\eqdef(v,\mapS(v))$ is bijective and $\hmapS(\sol)=\solh$.\end{lemma}
  \begin{sloppypar}
    The proof consists in observing that the convex function $\reg$ must coincide with the concave function $\minval-\fidp$ on $\sol$. We omit it for brevity. As a consequence of \Cref{faces:sec:linear}, we obtain ${j=\dim \face{0}{\sol} = \dim \face{\vecepi{0}{0}}{\solh}=\dim(\hS)}$.
  \end{sloppypar}

As for $\epir$ and $\lst{\reg}{0}$, we have $\epir\cap (\vecgal\times \{0\}) = \lst{\reg}{0}\times \{0\}$, and using~\eqref{eq:faceintercvx} and~\eqref{eq:faceprodcvx} we obtain 
\begin{align}
  \face{\vecepi{0}{0}}{\epir}\cap \left(\vecgal\times \{0\} \right) &=\face{\vecepi{0}{0}}{\epir\cap (\vecgal\times \{0\})}  = \face{0}{\lst{\reg}{0}}\times \{0\}.
\end{align}
Let $\Evec\eqdef\Span\left(\face{0}{\lst{\reg}{0}}\right)\subseteq \vecgal$. From \eqref{eq:spanintercvx}, we note that the linear spans $\hE$ and $\Evec$ are related through $\hE\cap \left(\vecgal\times \{0\} \right)= \Evec\times \{0\}$.
As a result 
\begin{align}
  \dim \hE&= \dim \Evec  +\codim_{\hE}(E\times\{0\}),
\end{align}
where $\codim_{\hE}(E\times\{0\})\in \{0,1\}$. If $\codim_{\hE}(E\times\{0\})=0$, we say that the face $\face{\vecepi{0}{0}}{\epir}$ is \emph{horizontal}. Otherwise we say that it is \emph{oblique}.

Now, we examine $\hypof$ and $\lst{\fido}{\fido(0)}$.  Let $Z$ be a linear complement to $\ker\Phi$ in $\vecgal$. Since $\rank\Phi =\nbmes$, the restriction $\restr{\Phi}{Z}:Z\rightarrow \RR^\nbmes$ is an isomorphism. As a result, the mapping
\begin{align*}
  \psi:  \begin{array}{rccrcl}
    \vecgal\times \RR & \longrightarrow & \ker\Phi\times Z\times \RR &\longrightarrow&(\ker\Phi)\times (\RR^\nbmes\times \RR)\\
    (u,r)&\longmapsto & (k,z,r)& \longmapsto&(k,(\Phi z,r))
 \end{array}
\end{align*}
(where $(k,z)$ is the unique element in $\ker\Phi\times Z$ such that $u=k+z$) is an isomorphism.
In particular, since $\psi$ maps $\hypof$ to $\ker\Phi\times (\hypo(\minval-\fido))$,
\begin{align*}
  \psi\left(\face{\vecepi{u}{r}}{\hypof}\right)&=\face{\psi\vecepi{u}{r}}{\psi(\hypof)}\\
  &=\face{(k,(\Phi z,r))}{\ker\Phi\times (\hypo(\minval-\fido)}\\
  &=\ker\Phi\times \face{(\Phi z,r)}{\hypo(\minval-\fido)}.
\end{align*}
Applying this to $(u,r)=(p,\reg(p))=(0,0)$ and considering the linear spans, we obtain
\begin{align}\label{faces:eq:codimH}
  \codim_{\vecgal\times\RR} \hH= \nbmes + 1 -\dimfidepi,
\end{align}
where $\dimfidepi$ is the dimension of $\face{\vecepi{0}{0}}{\hypo(\minval-\fido)}$.
Just like the faces of the regularizer, we may define $\Hvec=\Span \face{0}{\lst{\fido}{\fido(0)}}$ and notice that $\dim\hH=\dim\Hvec+\codim_{\hH}(\Hvec\times\{0\})$.
As a result, $\dimfidepi=\dim\face{0}{\lst{\fido}{\fido(0)}}$ (in which case we say the face $\face{\vecepi{0}{0}}{\hypo(\minval-\fido)}$ is \emph{horizontal}), or $\dimfidepi=\dim\face{0}{\lst{\fido}{\fido(0)}}+1$ (in which case we say it is \emph{oblique}).

From the above discussion, we see that the horizontality or obliqueness of the faces in $\epir$ and $\hypof$ play an important role. The following condition will be useful.
\begin{definition}\label{faces:def:oblique}
  We say that $p$ satisfies the \emph{double obliqueness condition} if both $\face{\vecreg{p}}{\epir}$ and $\face{\vecepi{\Phi{p}}{\minval-\fido(\Phi p)}}{\hypo(\minval-\fido)}$ are oblique. In other words, 
  \begin{align}
    \label{eq:oblique}
    \codim_{\hE}(\Evec\times \{\reg(p)\}) = 1 \qandq \codim_{\hH}(\Hvec\times \{\reg(p)\}) = 1.
  \end{align}
\end{definition}

\subsection{Conclusion of the proof of \texorpdfstring{\Cref{thm:first}}{Theorem \ref{thm:first}}}
It remains to study the last term of~\eqref{faces:eq:dimensions}.
\begin{lemma}\label{faces:lem:codim}
  The following inequality holds:  $\codim_{\vecgal\times\RR}(\hE+\hH)\geq 1$. 
  
  If both $\face{\vecepi{0}{0}}{\epir}$ and $\face{\vecepi{0}{0}}{\hypo(\minval-\fido)}$ are horizontal and ($\reg(0)> \inf \reg$ or $\fidp(0)>\inf \fidp$), then $\codim_{\vecgal\times\RR}(\hE+\hH)\geq 2$.
\end{lemma}

\begin{proof}
  First, in the general case, we prove that $\codim_{\vecgal}(\hE+\hH)\geq 1$ by arguing that $(\{0\}\times \RR)\cap (\hE+\hH)=\{0\}$. Indeed, assume by contradiction that $(\{0\}\times \RR)\subseteq (\hE+\hH)$. Then, for all $\epsilon>0$, there exists a vector $\vecepi{e}{\alpha}\in \hE$, a vector $\vecepi{h}{\beta}\in \hH$, such that $\vecepi{0}{-\epsilon} = \vecepi{e}{\alpha}+\vecepi{h}{\beta}$. Possibly reducing $\epsilon$, we may assume that $\vecepi{e}{\alpha}$ and $\vecepi{h}{\beta}$ are so small that $\vecepi{e}{\alpha}\in \face{\vecepi{0}{0}}{\epir}$ and  $\vecepi{h}{\beta}\in\face{\vecepi{0}{0}}{\hypof}$.
  Arguing as in \Cref{faces:sec:linear}, it is possible to prove that $\reg$ (resp.\ $(\minval-\fidp)$ coincides with a linear function on $\enscond{e'}{\vecepi{e'}{\alpha'}\in \face{\vecepi{0}{0}}{\epir}}$ (resp.\ on $\enscond{h'}{\vecepi{h'}{\beta'}\in \face{\vecepi{0}{0}}{\hypof}}$), with $\reg(e')=\alpha'$ and $\beta'=\minval-\fidp(h')$.

Hence,
  \begin{align*}
    (\reg+\fidp)(e)=\reg(e)+\fidp(-h)= \reg(e)+ (\fidp(-h)-\minval)+\minval=\alpha+\beta+\minval =\minval-\epsilon<\minval,
  \end{align*}
  which contradicts the fact that $\minval$ is the minimal value of \eqref{faces:eq:minconv}.

  Now, we assume that both $\face{\vecepi{0}{0}}{\epir}$ and $\face{\vecepi{0}{0}}{\hypof}$ are horizontal and  that $\reg(0)> \inf \reg$ (the case $\fidp(0)>\inf \fidp$ being similar), and we prove $\codim_{\vecgal\times\RR}(\hE+\hH)\geq 2$.
  Note that, both faces being horizontal, we have $\hE+\hH\subseteq (\vecgal\times \{0\})$. As a result $\hE=\Evec\times \{0\}$ and $\hH=\Hvec\times\{0\}$, and it is sufficient to prove that $\Evec+\Hvec$ has a nontrivial complement in $\vecgal$.
  By assumption, there exists $x'\in\vecgal$ such that $\reg(x')<\reg(0)$. By contradiction, if $\Evec+\Hvec=\vecgal$, we may write $x'=e+h$ with $e\in \Evec$, $h\in \Hvec$. For $\theta>0$, we consider the point 
  \begin{align*}
    x_{\theta}\eqdef \frac{\theta}{1+\theta}x'+ \frac{1}{1+\theta}(-\theta e) = \frac{\theta}{1+\theta}h. 
  \end{align*}
  Since $0$ is internal to $\face{0}{\lst{\reg}{0}}$ (resp. $\face{0}{\lst{\fidp}{\fidp(0)}}$, for $\theta>0$ small enough, $(-\theta e)\in \face{0}{\lst{\reg}{0}}$ (resp. $\frac{\theta}{1+\theta}h\in \face{0}{\lst{\fidp}{\fidp(0)}}$). Hence,
  \begin{align}
    (\reg+\fidp)(x_\theta)\leq \frac{\theta}{1+\theta}\underbrace{\reg(x')}_{<\reg(0)}+\frac{1}{1+\theta}\underbrace{\reg(-\theta e)}_{=\reg(0)}+\fidp(x_\theta)< (\reg+\fidp)(0),
     \end{align}
     which contradicts the fact that $0$ is a minimizer.\\
  As a result, $\Evec+\Hvec\subsetneq \vecgal$. If $W$ is a linear complement to $\Evec+\Hvec$ in $\vecgal$, then $(W\times \{0\})\oplus (\{0\}\times \RR)$ is thus a linear complement to $\hE+\hH$ in $\vecgal\times \RR$, with dimension $\geq 2$.
\end{proof}

  Now, we can finish the proof of \cref{thm:first}.
   From \eqref{faces:eq:dimensions} and \eqref{faces:eq:codimH}, we note that
\begin{align*}
  \dim \hE &= \nbmes +1-\dimfidepi +j -\codim_{\vecgal\times\RR}(\hE+\hH).
\end{align*}

    First, we assume that both $\face{\vecepi{0}{0}}{\epir}$ and $\face{\vecepi{0}{0}}{\hypo(\minval-\fido)}$ are horizontal: $\dim \hE =\dim \Evec$ and $\dimfidepi=\dimfidlst$. 
If $\reg(0)=\inf \reg$ and $\fidp(0)=\inf \fidp$, the inequality of \cref{faces:lem:codim} yields $\dim \Evec\leq \nbmes-\dimfidlst+j$. On the other hand, if $\reg(0)>\inf \reg$ or $\fidp(0)>\inf \fidp$, the stronger inequality of \cref{faces:lem:codim} yields $\dim\Evec \leq \nbmes-\dimfidlst+ j-1$

    \begin{sloppypar}    Now, if $\face{\vecepi{0}{0}}{\epir}$ is oblique, $\dim\hE=\dim \Evec+1$ and $\reg(0)>\inf \reg$. If $\face{\vecepi{0}{0}}{\hypo(\minval-\fido)}$ is oblique, then $\dimfidepi=\dimfidlst+1$ and $\fido(0)>\inf \fido$. Using that $\codim_{\vecgal}(\hE+\hH)\geq 1$, we obtain $\dim \Evec\leq \nbmes-\ell+j-1$ if only one of the faces is oblique, and $\dim \Evec\leq \nbmes-\ell+j-2$ if both are oblique. This concludes the proof.\end{sloppypar}

    \subsection{Proof of \texorpdfstring{\Cref{coro:repres}}{Corollary \ref{coro:repres}}}
  We derive the representation of $p=0$ as a convex combination. Let $\dimreg<+\infty$ be the above-mentioned upper-bound on $\dim \Evec$ and let $F$ be the linear closure of $\face{0}{\lst{\reg}{\reg(0)}}$. As $\face{0}{\lst{\reg}{\reg(0)}}\subseteq \Evec$, we know that $F$ is in fact its closure for the finite-dimensional topology of $\Evec$. By the Carathéodory-Klee theorem~\cite[Th. 3]{klee_theorem_1963}, since $F$ is closed and contains no line, its points are convex combinations of at most $\dimreg+1$ extreme points of $F$ or $\dimreg$ points, each an extreme point or a point in an extreme ray. To conclude, we note that $F$ is a face of $\lst{\reg}{\reg(0)}$, its extreme points are thus extreme points of $\lst{\reg}{\reg(0)}$.
  \qed 

\section{Level sets containing lines}
 \subsection{The generalized splines problem.} In~\cite{fisher_spline_1975}, Fisher and Jerome have in fact considered the more general problem
\begin{equation}\label{eq::spline}
  \min_{u\in \vecgal} \tvmes{\Ld u}  \quad \mbox{s.t.}\quad  \Ld u\in \Mm(\Omega) \qandq \Phi u =y,
\end{equation}
where $\vecgal\subseteq \mathcal{D}'(\Omega)$ is a suitably defined space of distributions, and the differential operator $\Ld:\mathcal{D}'(\Omega)\rightarrow \mathcal{D}'(\Omega)$ maps $\vecgal$ onto $\Mm(\Omega)$. We refer to~\cite{fisher_spline_1975} for precise assumptions and to \cite{unser2017splines,flinth_exact_2019,gupta_continuous-domain_2018} for extensions to a more general setting.
Under some topological assumptions, the results in~\cite{unser2017splines,flinth_exact_2019} describe the extreme points of the solution set to~\eqref{eq::spline} as generalized splines 
 \begin{align}\label{eq:fisherjeromekernel}
   u(x)&= \sum_{k=1}^{r} a_k\Ld^+\delta_{x_k}+ P(x)
\end{align}
where $a_k\in \RR$, $\Ld^+$ is some pseudo-inverse of $\Ld$ (one may see $\Ld^+\delta_{x_k}$ as a Green function for $\Ld$) and $P\in \Ker \Ld$. Moreover, $r\leq \nbmes -\dim \Phi(\ker\Ld)$.

For $\Omega=\RR$ and $\Ld=\mathrm{D}^n$, where $\mathrm{D}$ denotes the differentiation operator, one obtains the polynomial splines \begin{align*}
  u(x)&= \sum_{k=1}^{r} \frac{(x-x_k)_+^n}{n!} + P(x),
\end{align*}
where $P\in \RR_{n-1}[X]$.
In this section, we examine how our discussion can be extended to the case of level set containing lines so as to obtain representations of the form~\eqref{eq:fisherjeromekernel}.

  \subsection{Convex sets and their lineality space.} We first need to recall several properties of convex sets containing lines (see for instance~\cite{klee_extremal_1957} or \cite[Ch.8]{rockafellar_convex_1997}).
   We say that a nonempty convex set $C\subseteq \vecgal$ is invariant in the direction $v\in \vecgal$ if 
\begin{align}
  \label{eq:containline}
 \cvx+\RR v \subseteq \cvx.
\end{align} The collection of all vectors $v\in \vecgal$ such that~\eqref{eq:containline} holds is a vector space called the \emph{lineality space} of $\cvx$, denoted by $\lin{\cvx}$. 

If $\cvx$ is internal or linearly closed, given  $v\in \vecgal\setminus \{0\}$, it is equivalent to say that $C$ is invariant in the direction $v$, or to say that $C$ contains a line directed by $v$, \ie{} $(x_0+ \RR v)\subseteq C$ for some $x_0\in \vecgal$.
As a consequence, if $C_1$, $C_2$ are two nonempty convex sets, then $\lin{C_1}\cap \lin{C_2}\subseteq \lin{C_1\cap C_2}$, with equality if $C_1$ and $C_2$ are internal or linearly closed.


If $\cvx$ is invariant by some subspace $\linreg$, \ie{} $\linreg\subseteq\lin{\cvx}$, it is sometimes convenient to quotient the ambient space by $\linreg$. Considering a linear complement\footnote{We use freely the axiom of choice, hence any subspace of $\vecgal$ admits a complement subspace.} $W$ to $\linreg$, we note that there is a linear isomorphism 
\begin{align*}
  \psi:  \begin{array}{rccrcl}
    \vecgal & \longrightarrow & \linreg\times W &\longrightarrow& \linreg\times \left(\faktorf{\vecgal}{\linreg}\right)\\
    u&\longmapsto & (k,w)& \longmapsto&(k,\projreg(w))
 \end{array}
\end{align*}
where $(k,w)$ is the unique element in $ \linreg\times W $ such that $u=k+w$, and $\projreg:\vecgal\rightarrow \faktorf{\vecgal}{\linreg}$ is the canonical surjection. In other words $\projreg(w)=\projreg(u)\eqdef u+\linreg$ is the coset of $u$.

From~\eqref{eq:facebijcvx}, we note that $\psi\left(\face{u}{\cvx}\right)=\linreg\times \face{\projreg(u)}{\projreg(\cvx)}
$ so that the faces of $\projreg(\cvx)$ can be easily deduced from those of $\cvx$ and conversely, with 
\begin{align}\label{eq:facemoduloK}
  \projreg(\face{u}{\cvx})&=\face{\projreg(u)}{\projreg(\cvx)}.
\end{align}

If $\cvx$ is internal (resp. linearly closed), then $\projreg(\cvx)$ is internal (resp. linearly closed) and, if $\linreg=\lin{\cvx}$, $\projreg(\cvx)$ contains no line.


\subsection{Back to the optimization problem.}
Let $\linreg =\lin{\lst{\reg}{\reg(p)}}$ and $\linfid\eqdef \ker \Phi$. We note that \todetail{$\face{p}{\sol}$ is invariant by $\linreg\cap \linfid$}.
\begin{detail}
  Indeed, the face of the epigraph $\face{\vecreg{p}}{\epir}$ is internal and contains $\linreg\times\{\reg(p)\}$, hence it is invariant by $\hlinreg\eqdef \linreg\times\{0\}$. On the other hand, the hypograph $\hypof$ (hence $\face{\vecreg{p}}{\hypof}$) is invariant by $\hlinfid\eqdef \linfid\times \{0\}$. 
  From~\eqref{eq:faces} we deduce that $\face{\vecreg{p}}{\epir\cap \hypof}$ is invariant by $\hlinreg\cap \hlinfid = (\linreg\cap \linfid)\times \{0\}$. Since $\face{p}{\sol}$ is its projection onto the horizontal hyperplane (see \Cref{faces:sec:linear}),  it is invariant by $\linreg\cap \linfid$.
\end{detail}
Therefore, $\face{p}{\sol}$ is linearly isomorphic to ${(\linreg\cap \linfid)\times  \projregfid\left(\face{p}{\sol}\right)}$. Instead of considering the dimension of $\face{p}{\sol}$ to describe the point $p$, the following theorem relies on the dimension of the coset $\projregfid\left(\face{p}{\sol}\right)$.

\begin{theorem}
\label{thm:second}
  Let $\reg:\vecgal \to \RR\cup\{+\infty\}$, $\fido:\RR^\nbmes \rightarrow \RR\cup \{+\infty\}$ be two convex functions, and let $\Phi:\vecgal \rightarrow \RR^\nbmes$ be linear surjective. Assume that $p \in \sol$, with $\reg(p)+\fido(\Phi p)<+\infty$, and that $\lst{\reg}{\reg(p)}$ is linearly closed. Let $\linreg\eqdef \lin{\lst{\reg}{\reg(p)}}$, $d\eqdef \dim \Phi(\linreg)$, and $\linfid\eqdef \ker \Phi$.
  
  If $\dim\left( \projregfid\left(\face{p}{\sol}\right)\right)=j<+\infty$, then $\projreg(p)$ belongs to a face of $\projreg(\lst{\reg}{\reg(p)})$ with dimension at most $\dimreg$, where
  \begin{align}\label{faces:eq:dimregquotient}
    k\eqdef \begin{cases}
      \nbmes -\dimfidlst+j-d-1 & \mbox{if $\left(\reg(p)>\inf \reg\right)$ or $\left(\fido(\Phi p)>\inf \fido\right)$},\\
      \nbmes -\dimfidlst+j-d & \mbox{otherwise,}
    \end{cases}
  \end{align}
  and $\ell$ is the dimension of the minimal face of $\Phi p$ in $\lst{\fido}{\fido(\Phi p)}$.
  
  In particular, $\projreg(p)$ can be written as a convex combination of (at most):
\begin{itemize}
  \item[$\circ$] $\dimreg+1$ extreme points of $\projreg(\lst{\reg}{\reg(p)})$,
  \item[$\circ$] or $\dimreg$ points of $\projreg(\lst{\reg}{\reg(p)})$, each an extreme point or a point in an extreme ray.
\end{itemize}

  If, moreover, $p$ satisfies the obliqueness condition described in~\cref{faces:def:oblique}, the number $k$ can be reduced to $\nbmes -\dimfidlst+j-d-2$. 
\end{theorem}

 In particular, if  $p_1, \ldots, p_r\in \lst{\reg}{\reg(p)}$ are such that $\projreg(p_1),\ldots, \projreg(p_r) $ denote those extreme points (or points in extreme rays), 
\begin{equation}\label{eq:convcomb}
  p= \sum_{i=1}^r \theta_i p_i + u_{\linreg}, \qwhereq  \theta_i\geq 0,\  \sum_{i=1}^r \theta_i =1, \qandq u_{\linreg}\in \linreg.
\end{equation}

\begin{rem}
	In practice, if $\epir$ is linearly closed (\eg{} if $\reg$ is lower semi-continuous for some \revision{vector} topology), then the whole solution set $\sol$ is invariant by $(\linreg\cap \linfid)$, and $\projregfid\left(\face{p}{\sol}\right)=\left(\face{\projregfid(p)}{\projregfid(\sol)}\right)$. However, notice that the solution set $\sol$ may have more invariant directions than just $(\linreg\cap \linfid)$.
\end{rem}

\begin{proof}
  We assume, up to a change of the origin in $\vecgal\times\RR$, that $\vecreg{p}=\vecepi{0}{0}$. 

  As noted above, $\face{\vecreg{p}}{\epir}$ is invariant by  $\hlinreg=\linreg\times\{0\}$, hence $\hlinreg\subseteq \hE$. Similarly, $\hlinfid\subseteq \hH$, and classical results on quotient spaces~\cite[Ch.~3, Sec.~1]{lang_algebra_2002} imply  
\begin{align}
  &\left\{\begin{array}{cl}
    \faktor{\hE}{\hlinreg}&\approx\faktorf{\left(\faktor{\hE}{(\hlinreg\cap \hlinfid)}\right)}{\left(\faktor{\hlinreg}{(\hlinreg\cap\hlinfid)}\right)}\\
    \faktor{\hE}{(\hE\cap\hH)}&\approx\faktorf{\left(\faktor{\hE}{(\hlinreg\cap \hlinfid)}\right)}{\left(\faktor{(\hE\cap \hH)}{(\hlinreg\cap\hlinfid)}\right)}
  \end{array}\right.\\
  \intertext{hence, provided the corresponding dimensions are finite (we prove below that they are),}
  \dim(\faktorf{\hE}{\hlinreg})&=\dim\left(\faktorf{\hE}{(\hE\cap\hH)}\right) + \dim\left(\faktorf{(\hE\cap \hH)}{(\hlinreg\cap\hlinfid)}\right)-\dim\left(\faktorf{\hlinreg}{(\hlinreg\cap\hlinfid)}\right).
\end{align}

  The term $\dim\left(\faktorf{\hE}{(\hE\cap\hH)}\right)=\codim_{\hE}(\hE\cap\hH)$ has already been studied in~\Cref{sec:epireform,sec:facepi}. 
  Now, we prove that $\dim\left(\faktorf{(\hE\cap \hH)}{(\hlinreg\cap\hlinfid)}\right)=j$. Since $\hmapS$ (see \cref{faces:sec:linear}) satisfies $\hmapS(\linreg\cap \linfid)=(\linreg\cap \linfid)\times\{0\}$, it induces a map $\htmapS$ which makes the following \revision{commutative diagram}.
\begin{center}
  \begin{tikzcd}
    \Svec \arrow[r, "\hmapS"] \arrow[d, "\projregfid"]& \hS\arrow[d, "\hprojregfid"]\\
    \faktor{\Svec}{\linreg\cap \linfid} \arrow[r, "\htmapS"]& \faktor{\hS}{\hlinreg\cap\hlinfid}
  \end{tikzcd}
\end{center}
Since $\hmapS$ is bijective, $\htmapS$ is surjective by construction. Let us observe that \todetail{it is also injective}. 
\begin{detail}
  Consider any coset in $\faktorf{\Svec}{\linreg\cap \linfid}$, say $\projregfid(s)$ for some $s\in \Svec$. If $\htmapS(\projregfid(s))=0$, it means that $\hprojregfid(\hmapS(s))=0$, that is $(s,\mapS(s))\in (\linreg\times\linfid)\times\{0\}$. Hence $\projregfid(s)=0$, and $\htmapS$ is injective.
\end{detail}
  As a result $\dim\left(\faktorf{\hS}{\hlinreg\cap\hlinfid}\right)=\dim\left(\faktorf{\Svec}{\linreg\cap \linfid}\right)$. But by linearity, 
  \begin{align*}
\faktor{\Svec}{\linreg\cap \linfid}=\projregfid\left(\Svec\right)=\projregfid\left(\Span\left(\face{0}{\sol}\right)\right)=\Span\left(\projregfid\left(\face{0}{\sol}\right)\right), 
  \end{align*}
  so that $\dim\left(\faktorf{\hS}{\hlinreg\cap\hlinfid}\right)=j$.
  
  It remains to identify $\dim\left(\faktorf{\hlinreg}{(\hlinreg\cap\hlinfid)}\right)$. 
By the first isomorphism theorem, if $\restr{\Phi}{\linreg}$ denotes the restriction of $\Phi$ to $\linreg$,
\begin{align*}
  \Im \restr{\Phi}{\linreg}\approx \left(\faktor{\linreg}{\ker \restr{\Phi}{\linreg}}\right) = \left(\faktor{\linreg}{\linreg\cap \linfid}\right)\approx \left(\faktor{\hlinreg}{\hlinreg\cap \hlinfid}\right). 
\end{align*}
  As a result, $\dim (\faktorf{\hlinreg}{\hlinreg\cap \hlinfid})=\dim \Phi(\linreg)=d$.

To conclude, we deal with $\hE$. As  in \Cref{sec:facepi}, we note that $\dim(\faktorf{\hE}{\hlinreg}) = \dim(\faktorf{\Evec}{\linreg})+\mbox{($0$ or $1$)}$, depending on whether $\face{0}{\epir}$ is horizontal or oblique. 
 Hence,  $\dim(\faktorf{\Evec}{\linreg})\leq \dimreg$ where $\dimreg$ is given by \eqref{faces:eq:dimregquotient}.
  Moreover, 
  \begin{align*}
    \faktor{\Evec}{\linreg}= \projreg(\Span(\face{0}{\lst{\reg}{0}})) = \Span\left(\projreg(\face{0}{\lst{\reg}{0}})\right)
  \end{align*}
  so that $\dim(\faktorf{\Evec}{\linreg})$ is the dimension of the elementary face of $\projreg(p)$ in $\projreg(\lst{\reg}{\reg(p)})$. As $\lst{\reg}{\reg(p)}$ is linearly closed and $\linreg$ is its lineality space, we note that $\projreg(\lst{\reg}{\reg(p)})$ is linearly closed and contains no line. The representation of $\projreg(p)$ as a convex combination follows from the same argument as in the proof of \cref{coro:repres}.  
  
\end{proof}

\subsection{Link with the Fisher-Jerome result.}
The conclusions of~\Cref{thm:second} and \eqref{eq:convcomb} recover the representation result~\eqref{eq:fisherjeromekernel} for the Fisher-Jerome problem. Indeed, assuming that there is a point $p$ such that $j=0$ (\eg{} if $\projregfid(\sol)$ is compact for some \revision{locally convex Hausdorff} topology), we obtain that $\projreg(p)$ belongs to a face with dimension at most $\nbmes-d-1$. The quantities $\frac{a_k}{\abs{a_k}}\Ld^+\delta_{x_k}$ are some specific choices of  points $p_k$ such that $\projreg(p_k)=\frac{a_k}{\abs{a_k}}\delta_{x_k}$ is an extreme point of the total variation unit ball. The operator $\Ld^+$ is determined by choosing some specific linear complement $W$ to $\linreg$. Conversely, if $W$ (hence $\Ld^+$) is fixed and $p_1, \ldots, p_r$ are given as in~\eqref{eq:convcomb}, one can recover~\eqref{eq:fisherjeromekernel} by adding some elements of $\linreg$ (\ie{} changing $u_\linreg$).

\section{Is it optimal?}
\subsection{Optimality in the general class of convex problems.} 
For each different case of \Cref{thm:first} and \Cref{coro:repres}, it is possible to exhibit a convex function for which the upper-bounds are attained. In that sense, the results stated above are optimal. Moreover, due to their geometric essence and the weakness of their \revision{topological} assumptions\footnote{\revision{The notions of linearly closed set and algebraic interior used in the present paper can be regarded as topological, since a set is linearly closed (resp. algebraically open) if and only if it is closed (resp. open) in the \emph{core topology} (see \cite{klee_jr.1951}). However, in practice, that topology is rarely used since it is not compatible with the vector space structure of $\vecgal$, so that one might as well regard those notions as algebraic definitions. In some cases, a convex set is linearly closed if and only if it is closed for the \emph{convex core topology}, \ie the finest locally convex topology (see~\cite[Prop. 4.6 and 8.7]{klee_jr.1951}). In any case, linear closedness is a weak assumption, since it is implied by being closed in a topological vector space. On the other hand, a convex set is algebraically open if and only if it is open for the convex core topology.}}, they are quite general. In comparison, an approach relying on the subdifferential and optimality conditions would require a constraint qualification argument together with a suitable choice of topology. Moreover, it is not clear to the author whether exploiting the subdifferential could describe precisely the extreme points that are not exposed. 

On the other hand, for several specific problems that we discuss below, the bounds provided by \Cref{thm:first} are not tight.

\subsection{The Carath\'eodory number.} A first example, discussed in~\cite{boyer_representer_2019}, is the case of Semi-Definite Programs (SDP), such as the feasibility problem
\begin{align*}
  \min_Q \chi_{\psdr{n}}(Q)  \quad\mbox{s.t.}\quad \Phi Q=y,
\end{align*}
where $\psdr{n}$ denotes the cone of symmetric positive semi-definite matrices of size $n$ \revision{and, for a given set $A$, $\chi_A(x)\eqdef 0$ if $x \in A$, and $+\infty$ otherwise}.
While \Cref{coro:repres} predicts that an extreme point of the solution set is a convex combination of at most $\nbmes$ rank-one matrices (\ie{} the points in the extreme rays of $\psdr{n}$), it is known~\cite{barvinok_problems_1995} that there is a solution which is made of at most $\frac{1}{2}\left(\sqrt{8\nbmes+1}-1\right)\leq \nbmes$ such atoms.
As it turns out, even though the solution does belong to a face of $\psdr{n}$ with dimension $\nbmes$, the Carath\'eodory-Klee theorem is not sharp, since \emph{one needs less than $\nbmes$ points in extreme rays to represent the points of that face.}
That phenomenon is quite common with non-polyhedral sets: for instance, in the Euclidean closed unit ball, every point is a convex combination of at most two extreme points (regardless of the ambient dimension). More generally, if $F\subseteq \vecgal$ is convex, linearly bounded, linearly closed and each point of $F$ is internal or an extreme point (\eg{} if $F$ is strictly convex), the same property holds.

On the other hand, if $F\subseteq \vecgal$ is an $\nbmes$-dimensional convex polytope, it is possible to check that almost every  point of $F$ (in the sense of the Lebesgue measure) is a convex combination of $\nbmes+1$ (and not less) extreme points of $F$.

\subsection{Nonnegative measures.} Another interesting case is the truncated trigonometric moment problem,
\begin{align}\label{eq:caratoeplitz}
  \min_\mes \chi_{\Mm^+(\TT)}(\mes)  \quad\mbox{s.t.}\quad \int_\TT \varphi_k(t)\D\mes(t)=y_k\quad (0\leq k\leq 2f_c),
\end{align}
where $\Mm^+(\TT)$ is the set of nonnegative measures on the torus $\TT=\RR/\ZZ$, $y\in \RR^{2f_c+1}$, and 
$\varphi_0(t)=1$, $\varphi_{2j-1}(t)= \cos(j 2\pi t)$ and $\varphi_{2j}(t)= \sin(j 2\pi t)$ for $1\leq j\leq f_c$.

The Carath\'eodory-Toeplitz theorem \cite{caratheodory_uber_1907,toeplitz_uber_1911} states that there is a solution to~\eqref{eq:caratoeplitz} if and only if the matrix
\begin{align*}
  T(c)\eqdef  \begin{pmatrix}
    c_0 & c_1 & \cdots & c_{f_c}\\
    c_1^* &c_0 & \ddots & \vdots\\
    \vdots & \ddots&  \ddots&c_1 \\
    c_{f_c}^* &\cdots &c_1^* & c_0 
  \end{pmatrix},\qwhereq  c_j\eqdef y_{2j-1}-\imath y_{2j},\ c_0\eqdef y_0,
\end{align*}
is positive semi-definite. If $r\eqdef \rank T(c)\leq f_c$, the solution $\mes$ is unique, and its support has cardinality $r$. If $T(c)$ is invertible, there is an infinity of solutions with cardinality $f_c+1$ (and more). In particular for any $t_0\in \TT$ there is a solution which charges $\{t_0\}$. Note that similar results hold for T-systems on an interval~\cite[Ch. 4, Sec. 4]{krein_markov_1977}.

That result contrasts with \Cref{coro:repres} which would predict a sum of at most $2f_c+1$ Dirac masses. Here, the situation is different from the case of $\psdr{n}$, since any measure belonging to a $d$-dimensional elementary face of $\Mm^+(\TT)$ is a sum of \emph{exactly} $d$ Dirac masses. Therefore  we must have $\dim\face{m}{\Mm^+(\TT)}<2f_c+1$ and, recalling \eqref{faces:eq:dimensions}, we deduce that the lower bound on $\codim_{\vecgal\times \RR}(\hE+\hH)$ provided by \Cref{faces:lem:codim} is far too pessimistic.
In other words, the affine spaces determined by the Fourier coefficients only intersect very specific faces of the cone $\Mm^+(\TT)$. 

An intuitive explanation consists in counting the ``degrees of freedom'' of the problem (we do not consider the statistical notion used in~\cite{2019-Poon-DOF}, but simply the ``number of variables that should be fixed''). Informally, to fix the positions and amplitudes of $k$ Dirac masses, that is $2k$ variables, we need at least $2k$ equations, \ie{} $2k\leq 2f_c+1$.

\subsection{The Basis-Pursuit for measures.} Surprisingly, things are different when considering the \emph{Basis Pursuit for measures},
\begin{align}\label{eq:bptrigo}
  \min_{\mes\in \Mm(\TT)} \abs{\mes}(\TT)  \quad\mbox{s.t.}\quad \int_\TT \varphi_k(t)\D\mes(t)=y_k\quad (0\leq k\leq 2f_c),
\end{align}
where $\{\varphi_k\}_{k=0}^{2f_c}$ is again the trigonometric system. In~\cite{condat_atomic_2019}, Laurent Condat observed that when $y$ is the Fourier coefficient vector of two opposite close spikes, a solution to~\eqref{eq:bptrigo} is a Dirac comb. More precisely, if $y=\Phi \mes_0\eqdef \left(\int_{\TT}\varphi_k\D\mes_0\right)_{k=0}^{2f_c}$ with $\mes_0=\delta_{\h/2} -\delta_{-\h/2}$ and $0<h< \frac{1}{2f_c}$, he noticed that the solution to \eqref{eq:bptrigo} is given by
\begin{align}\label{eq:solblasso}
  \mes = \sum_{j=-f_c}^{f_c-1}a_j\delta_{t_j}, \qwhereq t_j\eqdef \frac{1}{4f_c}+\frac{j}{2f_c}, \ \mbox{and }\qquad\qquad\\
  a_j= (-1)^j \frac{\cos(\pi\h f_c)}{2f_c}\left(\cotan(\pi(\frac{1}{4f_c}+\frac{j}{2f_c}-\h/2)) - \cotan(\pi(\frac{1}{4f_c}+\frac{j}{2f_c}+\h/2))\right).\label{eq:amplblasso}
\end{align}
 While this observation in~\cite{condat_atomic_2019} seems to rely on numerical experiments, we provide in \Cref{sec:proofblasso} a proof relying on a duality argument. 

\begin{proposition}\label{prop:proofblasso}
  For $y=\Phi\mes_0$ with $\mes_0=\delta_{\h/2} -\delta_{-\h/2}$ and $0<\h<1/(2f_c)$, the unique solution to~\eqref{eq:bptrigo} is given by~\eqref{eq:solblasso} and~\eqref{eq:amplblasso}.
\end{proposition}

As a consequence of \Cref{prop:proofblasso}, the number of Dirac masses predicted by \Cref{coro:repres} is \emph{almost optimal} ($2f_c+1$ Dirac masses are predicted whereas $2f_c$ actually appear). In fact, one cannot do ``better'': it is proved in \cite{condat_atomic_2019} that for every $y\in \RR^{2f_c+1}$, there is a solution to~\eqref{eq:bptrigo} which is a sum of at most $2f_c$ Dirac masses. 

Arguing informally in terms of ``degrees of freedom'', we see that the above situation is quite peculiar: the relative positions of the Dirac masses are fixed, they can only move by a global translation. As a result, the $2f_c+1$ variables determine the $2f_c$ amplitudes of the spikes and the last degree of freedom which is a global shift of the Dirac comb.

\section{Conclusion}\label{conclusions}
The representer theorem presented in this note describes, under very weak assumptions, the solutions of variational problems as a combination of a few atoms. It is optimal in the sense that for each described configuration, there is a variational problem for which the predicted number of atoms is attained. However, in some specific problems, the solutions may be sparser than predicted by the theorem. It is for instance the case of the truncated trigonometric moment problem involving the positivity constraint in the space of measures. On the other hand, with the total variation regularization for signed Radon measures, the prediction of the theorem is almost optimal.

\section*{Acknowledgements}
The author warmly thanks Claire Boyer, Antonin Chambolle, Yohann De Castro, Frédéric de Gournay and Pierre Weiss for fruitful discussions about representer theorems. \revision{He also thanks an anonymous referee for stimulating comments about the connection between linear closedness and closedness in a topological sense.} 

\revision{This work was supported by the ANR CIPRESSI project, grant ANR-19-CE48-0017-01 of the French Agence Nationale de la Recherche.}

\appendix
\section{Proof of \texorpdfstring{\Cref{prop:proofblasso}}{Proposition \ref{prop:proofblasso}}}\label{sec:proofblasso}

We endow $\Mm(\TT)$ with the weak-* topology, and $\CT$ with the norm topology, so that $\Mm(\TT)$ and $\CT$ are \emph{paired spaces} (in the sense of, \eg{}~\cite{rockafellar_conjugate_1989}). A dual problem to~\eqref{eq:bptrigo} is 
\begin{align}
  \label{eq:dualbptrigo}
  \sup_{p\in \RR^{2f_c+1}} \int_{\TT}\left(\sum_k p_k \varphi_k(t)\right) \D m_0(t) \quad\mbox{s.t.}\quad \normi{\sum_k p_k \varphi_k}\leq 1.
\end{align}
It can be shown that strong duality holds (see~\cite{duval_exact_2015}) and that \eqref{eq:dualbptrigo} has a solution. Moreover, any admissible measure $\mes$ is a solution to \eqref{eq:bptrigo} and $p$ is a solution to~\eqref{eq:dualbptrigo} if and only if they satisfy the following extremality condition: the trigonometric polynomial $\eta\eqdef \sum_k p_k \varphi_k$ satisfies $\eta(t)=1$ for all $t\in \supp(\mes_+)$ and $\eta(t)=-1$ for all $t\in \supp(\mes_-)$, where $\mes_+$ and $\mes_-$ denote the positive and negative parts in the Jordan decomposition of $\mes$. Therefore, it suffices to solve \eqref{eq:dualbptrigo} to discover the possible support and sign of all the solutions to \eqref{eq:bptrigo}.

\begin{lemma}\label{lem:soldual}
  For $0<\h<\frac{1}{2f_c}$, the unique solution to \eqref{eq:dualbptrigo} is $p_*=(0,\ldots,0,1)$, corresponding to the function $\etas:t\mapsto \sin(2\pi f_c t)$.
\end{lemma}
\begin{proof}
  As $\mes_0=\delta_{\h/2} -\delta_{-\h/2}$, \Cref{eq:dualbptrigo} reformulates as the maximization problem
  \begin{align}\label{eq:maxeta}
  \max_{\eta} \left(\eta(\h/2) - \eta(-\h/2) \right)
\end{align}
where the maximization is over all the trigonometric polynomials $\eta$ with degree at most $f_c$ and ${\normi{\eta}\leq 1}$.
  Let $\eta$ be a solution to \eqref{eq:maxeta}. Note that its odd part, $\etao(t)\eqdef \frac{1}{2}\left(\eta(t)-\eta(-t)\right)$, is also a solution to $\eqref{eq:maxeta}$, so we first study $\etao$.

 Since $\normi{\etao}\leq 1$ and $\etas(t_j)=(-1)^{j}$ for $j\in \{-f_c,\ldots,f_c-1\}$, we note that $(-1)^{j}(\etas-\etao)(t_j)\geq 0$. Moreover,  
  \begin{align*}
    \etao(\h/2)=\frac{1}{2}\left[\etao(\h/2)-\etao(-\h/2) \right]\geq \frac{1}{2}\left[\etas(\h/2)-\etas(-\h/2)\right]=\etas(\h/2),
  \end{align*}
  and by oddness $\etao(-\h/2)\leq \etas(-\h/2)$. As a result, if we define 
  \begin{align*}
    \forall j\in \{-f_c-2\ldots, f_c-1\},\quad    t'_{j}\eqdef\begin{cases}
      t_{j+2} & \mbox{for $-f_c\leq j\leq -3$,}\\
      -\h/2 & \mbox{for $j=-2$,}\\
      \h/2 & \mbox{for $j=-1$,}\\
      t_j & \mbox{for $0\leq j\leq f_c-1$,}
    \end{cases}\\
  \end{align*}
  we have $2f_c+2$ distinct points such that $(-1)^{j}(\etas-\etao)(t'_j)\geq 0$. Thus, $\etas-\etao$ has at least $2f_c+2$ roots (counting multiplicity): it is clear by the mean value theorem if each inequality is strict, and it can also be checked by approximation and counting the roots (with multiplicity) of the limit in the general case. 
   As a consequence, $\etas-\etao$ has at least $2f_c+2$ roots and degree (at most) $f_c$: it must be identically zero.

We have proved that $\etao=\etas$, it remains to deal with the even part: we write $\eta=\etas+\etae$, where $\etae$ is the even part of $\eta$. Since $(-1)^j\eta(t_j)\leq 1$ for $-f_c\leq j\leq f_c-1$, we get
  \begin{align*}
    1+(-1)^j\etae(t_j)\leq 1,  
  \end{align*}
  so that $(-1)^j\etae(t_j)\leq 0$. Since $\etae(t_j)=\etae(t_{-j-1})$ we deduce that in fact $\etae(t_j)=0$, so that $\etae$ has $2f_c$ roots.
But, since $\eta$ and $\etas$ reach their maximum or minimum at each $t_j$, we also have $0=\eta'(t_j)=\etas'(t_j)+\etae'(t_j)=\etae'(t_j)$, so that $t_j$ is a double root of $\etae$. Hence $\etae=0$ and $\eta=\etas$.
  \end{proof}
As a consequence, the support of every solution to \eqref{eq:bptrigo} is contained in $\{t_{-f_c}, \ldots, t_{f_c-1}\}$.
Since that set is equispaced in $\TT$, we notice as in~\cite{condat_atomic_2019} that, to recover the amplitudes $a_j$, we may invert the system
\begin{align}\label{eq:fourierval}
  \sum_{j=-f_c}^{f_c-1} a_j e^{-2\imath \pi k t_j} &= c_k \quad (-f_c\leq k\leq f_c-1)
\end{align}
(where $c_0=y_0$, $c_k =y_{2k-1}-\imath y_{2k}$ for $k\geq 1$ and $c_{-k}=c_k^*$ otherwise) using the Discrete Fourier Transform.
We obtain
\begin{align*}
  a_\ell &= \frac{1}{2f_c}\sum_{k=-f_c}^{f_c-1} e^{2\imath \pi k\left(\frac{1}{4f_c}+\frac{\ell}{2f_c}\right)}c_k.
\end{align*}
Recalling that $c_k= e^{-2\imath \pi k\frac{\h}{2}} -e^{2\imath \pi k\frac{\h}{2}}$, we introduce the function
\begin{align}
  f(x)\eqdef \frac{1}{2f_c}\sum_{k=-f_c}^{f_c-1}e^{2\imath\pi kx} = e^{-\imath\pi x}\frac{\sin(2\pi f_c x)}{2f_c\sin(\pi x)} =\frac{\sin(2\pi f_c x)}{2f_c}\left(\cotan(\pi x)-\imath\right),
\end{align}
so that $a_\ell=f(\frac{1}{4f_c}+\frac{\ell}{2f_c}-\h/2)-f(\frac{1}{4f_c}+\frac{\ell}{2f_c}+\h/2)$.

Since $\sin\left(\frac{\pi}{2}+\pi \ell+\pi \h f_c \right) =\sin\left(\frac{\pi}{2}+\pi \ell-\pi \h f_c \right)=(-1)^\ell\cos(\pi\h f_c)$, we get
\begin{align*}
  a_\ell&= (-1)^\ell \frac{\cos(\pi\h f_c)}{2f_c}\left(\cotan(\pi(\frac{1}{4f_c}+\frac{\ell}{2f_c}-\h/2)) - \cotan(\pi(\frac{1}{4f_c}+\frac{\ell}{2f_c}+\h/2))\right).
\end{align*}
Since the function $\cotan$ is (strictly) decreasing on $\oi{0}{\pi}$, we see that $a_\ell\neq 0$ with $\sign(a_\ell)=\sign(\eta_*(t_\ell))$, which is the desired optimality condition.

It remains to check that \eqref{eq:fourierval} also holds for $k=f_c$. Since $a_\ell\in \RR$ for all $\ell$, we obtain it by taking the conjugate of \eqref{eq:fourierval} for $k=-f_c$.

To summarize, there is only one measure $\mes$ such that $\left(\int_{\TT}\varphi_k\D\mes(t)\right)_k=y$ and $\etas(t)=1$ for all $t\in \supp(\mes_+)$ and $\etas(t)=-1$ for all $t\in \supp(\mes_-)$. It is given by \eqref{eq:solblasso}. The extremality conditions imply that it is the only solution to \eqref{eq:bptrigo}.

\begin{rem}
  For $\h=\frac{1}{2f_c}$, a slight variation on the argument shows that the conclusions of~\Cref{lem:soldual} still hold, and that the unique solution  $\mes$ satisfies $a_0=1$, $a_{-1}=-1$ and $a_\ell=0$ otherwise. In other words, the measure $\mes_0= \delta_{\h/}-\delta_{-\h/2}$ is perfectly recovered by~\eqref{eq:bptrigo}.
\end{rem}

\bibliographystyle{alpha}
\bibliography{convex}

\end{document}

%% file: fig-epigraphs.tex
\tdplotsetmaincoords{85}{22}
\begin{tikzpicture}[scale=1.4, tdplot_main_coords,axis/.style={->,dashed},thick]

\tikzset{cross/.style={path picture={%
    \draw 
    (path picture bounding box.south east)--(path picture bounding box.north west) 
    (path picture bounding box.south west)--(path picture bounding box.north east); 
  }},cross/.default={0.5pt}} 


\def\h{2.5} 
\def\xo{0} 
\def\yo{0} 
\def\zo{-1} 

\def\k{5} 

\draw[thick,->,>=stealth] (\xo,\yo,\zo) -- (\xo+2.1*\h,\yo,\zo) node[anchor=north east]{};
\draw[thick,->,>=stealth] (\xo,\yo,\zo) -- (\xo,\yo+2.9*\h,\zo) node[anchor=north west]{};
\draw[thick,->,>=stealth] (\xo,\yo,\zo) -- (\xo,\yo,\zo+1+\h) node[anchor=south,left]{$\{0\}\times \RR$};

\coordinate  (k0) at (\xo,\yo,\zo){};
\coordinate  (k1) at (\xo+\k,\yo,\zo){};
\coordinate  (k2) at (\xo\k,\yo+\k,\zo){};
\coordinate  (k3) at (\xo,\yo+\k,\zo){};
\fill[gray!80,opacity=0.8] (k0) -- (k1) -- (k2) -- (k3) --  cycle; 
\coordinate  (k4) at (\xo+0.8*\k,\yo+0.5*\k,\zo){};
\draw (k4) node[] {$\vecgal\times\{0\}$};

\input{fig-epiR}

\input{fig-hypof}


\coordinate  (f0) at (\xshift+\a*\maxfidzface,\minregyface,\zo){}; 
\coordinate  (f1) at (\xshift+\a*\minfidzface,\minregyface,\zo){}; 
\coordinate  (g0) at (\xshift+\a*\maxfidzface,\maxregyface,\zo){}; 
\coordinate  (g1) at (\xshift+\a*\minfidzface,\maxregyface,\zo){}; 
\draw[regcolor!50!Bleu] (f0) -- (f1) -- (g1) -- (g0) -- cycle;
\fill[regcolor!50!Bleu,opacity=0.5] (f0) -- (f1) -- (g1) -- (g0) -- cycle;

\coordinate  (i0) at (\xshift+\a*\maxfidzface,\minregyface,\maxfidzface){}; 
\coordinate  (i1) at (\xshift+\a*\minfidzface,\minregyface,\minfidzface){}; 
\coordinate  (j0) at (\xshift+\a*\maxfidzface,\maxregyface,\maxfidzface){}; 
\coordinate  (j1) at (\xshift+\a*\minfidzface,\maxregyface,\minfidzface){}; 
\draw[regcolor!50!Bleu] (i0) -- (i1) -- (j1) -- (j0) -- cycle;
\fill[regcolor!50!Bleu,opacity=0.5] (i0) -- (i1) -- (j1) -- (j0) -- cycle;

  \node[minimum size=2pt] (P) at ($ (i0)!0.5!(j1) $) {\footnotesize$\bullet$};
  \draw[] (P) node[xshift=-16,yshift=12] {$\vecepi{p}{\reg(p)}$};

\draw[regcolor!50!Bleu,dotted,opacity=0.3] (f0) -- (i0);
\draw[regcolor!50!Bleu,dotted,opacity=0.3] (f1) -- (i1);
\draw[regcolor!50!Bleu,dotted,opacity=0.3] (g0) -- (j0);
\draw[regcolor!50!Bleu,dotted,opacity=0.3] (g1) -- (j1);
\draw[regcolor!50!Bleu] (g0) node[right] {$\sol$};
\draw[regcolor!50!Bleu] (j0) node[right] {$\solh$};

\node[minimum size=2pt] (PR) at ($ (f0)!0.5!(g1) $) {\footnotesize$\bullet$};
\draw[] (PR) node[yshift=-15] {$\vecepi{p}{0}$};
\draw[dotted,opacity=0.1] (P) -- (PR);

\end{tikzpicture}